\documentclass[notitlepage]{scrartcl}
\usepackage[shortlabels]{enumitem}
\usepackage{amsmath}
\usepackage{amssymb}
\usepackage{amsthm}
\usepackage{abstract}
\usepackage{xcolor}
\usepackage{comment}
\usepackage{mathtools}
\usepackage{graphicx}
\usepackage{caption}
\usepackage{subcaption}
\usepackage{float}
\usepackage{hyperref}
\usepackage{cleveref}
\hypersetup{
    colorlinks=true,
    linkcolor=blue,
    filecolor=magenta,      
    urlcolor=cyan
    }
\makeatletter
\newcommand*{\barfix}[2][.175ex]{%
  \mathpalette{\@barfix{#1}}{#2}%
}
\newcommand*{\@barfix}[3]{%
  \vbox{%
    \kern#1\relax
    \hbox{$#2#3\m@th$}%
  }%
}
\makeatother

\newtheorem{theorem}{Theorem}

\newtheorem{corollary}[theorem]{Corollary}
\newtheorem{lemma}[theorem]{Lemma}

\newtheorem{conjecture}[theorem]{Conjecture}

\newcommand{\footremember}[2]{%
    \footnote{#2}
    \newcounter{#1}
    \setcounter{#1}{\value{footnote}}%
}
\newcommand{\footrecall}[1]{%
    \footnotemark[\value{#1}]%
} 

\usepackage[margin=1in]{geometry} 

\title{\vspace{-1.5cm}Long cycles in percolated expanders}
\author{%
Maur\'{i}cio Collares \footremember{trailer}{\scriptsize{Institute of Discrete Mathematics, Graz University of Technology, Steyrergasse 30, 8010 Graz, Austria. Emails: mauricio@collares.org, erde@math.tugraz.at.}}
\and Sahar Diskin \footremember{alley}{\scriptsize{School of Mathematical Sciences, Tel Aviv University, Tel Aviv 6997801, Israel. Emails: sahardiskin@mail.tau.ac.il, krivelev@tauex.tau.ac.il.}}%
\and Joshua Erde \footrecall{trailer}%
\and Michael Krivelevich \footrecall{alley}%
}

\begin{document}
\maketitle

\vspace{-1cm}
\begin{abstract}
Given a graph $G$ and probability $p$, we form the random subgraph $G_p$ by retaining each edge of $G$ independently with probability $p$. Given $d\in\mathbb{N}$ and constants $0<c<1, \varepsilon>0$, we show that if every subset $S\subseteq V(G)$ of size exactly $\frac{c|V(G)|}{d}$ satisfies $|N(S)|\ge d|S|$ and $p=\frac{1+\varepsilon}{d}$, then the probability that $G_p$ does not contain a cycle of length $\Omega(\varepsilon^2c^2|V(G)|)$ is exponentially small in $|V(G)|$. As an intermediate step, we also show that given $k,d\in \mathbb{N}$ and a constant $\varepsilon>0$, if every subset $S\subseteq V(G)$ of size exactly $k$ satisfies $|N(S)|\ge kd$ and $p=\frac{1+\varepsilon}{d}$, then the probability that $G_p$ does not contain a path of length $\Omega(\varepsilon^2 kd)$ is exponentially small. We further discuss applications of these results to $K_{s,t}$-free graphs of maximal density.
\end{abstract}

\section{Introduction}
The binomial random graph $G(n,p)$ can be seen as a model of a random subgraph of the complete graph $K_n$, obtained by retaining each edge of $K_n$ independently with probability $p$. A particularly interesting feature of this model, first observed by Erd\H{o}s and R\'{e}nyi~\cite{ER59}, is the \textit{phase transition} that it undergoes with respect to its component structure when $p$ is around $\frac{1}{n}$. In~\cite{ER59} they showed that when $p \leq \frac{1-\varepsilon}{n}$, where $\varepsilon >0$ is some arbitrary constant, then with high probability (\textbf{whp}) every component of the graph is \emph{small}, of logarithmic order, whereas when $p \geq \frac{1+\varepsilon}{n}$, \textbf{whp} there is a unique \textit{giant} component of linear order. Later, Ajtai, Koml\'{o}s, and Szemer\'{e}di~\cite{AKS81a} showed that when $p \geq \frac{1+\varepsilon}{n}$, \textbf{whp} $G(n,p)$ contains a cycle of length $\Theta(n)$, settling a long-standing conjecture of Erd\H{o}s. For more background on the theory of random graphs, see~\cite{B01,FK16,JLR00}.

More recently, the emergence of similar structures in generalisations of this model has been studied. Given a host graph $G$, consider the random subgraph $G_p$ obtained by retaining each edge of $G$ independently with probability $p$. What (minimal) structural assumptions on the host graph $G$ are sufficient to guarantee the likely emergence of particular structures past certain ``natural'' thresholds?

Traditionally, these assumptions have been ``local'', in terms of the degree sequence of the graph, and in particular its minimum degree. Indeed, a standard coupling argument with a Galton--Watson branching process implies that for any graph $G$ with minimum degree $\delta(G) \geq d$ and $\varepsilon >0$, \textbf{whp} $G_p$ will contain a component of order $\Omega(d)$ when $p \geq \frac{1+\varepsilon}{d}$. The existence of long paths in this model and regime is less obvious. Through an analysis of the Depth-First Search algorithm, Krivelevich and Sudakov~\cite{KS13} showed that if $\delta(G)\ge d$ and $\varepsilon>0$, \textbf{whp} $G_p$ will contain a path of length $\Omega(d)$ when $p\ge \frac{1+\varepsilon}{d}$. However, in these general models, where the girth of $G$ could be much larger than the minimum degree, it is much less clear how to show the typical existence of a long \textit{cycle}. Using more delicate methods, this was eventually proven by Krivelevich and Samotij~\cite{KS14} (see also~\cite{R14}), who further conjectured that similar statements should hold when the bound on the minimum degree is replaced with a bound on the average degree.

Recent work of Diskin, Erde, Krivelevich, and Kang~\cite{DEKK23+} on the connection between the phase transition in $G_p$ and the isoperimetric properties of the host graph $G$ suggests another way of thinking about the quantitative aspects of the phase transition.
\begin{theorem}[{\cite[Theorem 4]{DEKK23+}}]\label{t:DEKK}
Let $k=\omega(1), d\leq k$ and let $G$ be a graph on more than $k$ vertices, such that every $S\subseteq V(G)$ with $|S|\leq k$ satisfies $e(S,V(G)\setminus S) \geq d|S|$. Let $\varepsilon > 0$ be a small constant and let $p = \frac{1+\varepsilon}{d}$. Then, with probability tending to $1$ as $k$ tends to infinity, $G_p$ contains a component of order at least $\frac{k}{2}$.
\end{theorem}
In fact, with slightly more assumptions on the graph, and a slightly weaker bound on the probability, a quantitatively similar statement to \Cref{t:DEKK} holds when the ``global'' assumption of expansion at all small scales is replaced by the assumption that sets of size \emph{exactly} $k$ expand well, see~\cite[Theorem 3.1]{DEKK23+}. 

One can think of Theorem \ref{t:DEKK} as giving an alternative heuristic for the nature of the phase transition. Indeed, here the point of criticality is controlled by the \textit{expansion ratio of subsets} (that is, $d$), whereas the quantitative aspects of the component structure above the critical point are controlled by the \emph{scale} at which this level of expansion holds (that is, $k$). 

Returning to the conjecture of Krivelevich and Samotij~\cite{KS13}, let us note that the assumption that sets of size exactly $k$ expand by a factor of $d$ is, in a sense, a strengthening of the assumption of having an average degree $d$ (indeed it implies the average degree is at least $d$). Nevertheless, for large $k$ this assumption is still somewhat weaker than a constraint on the minimum degree (in particular, it implies that there are at most $k$ vertices of degree less than $d$).
It is thus perhaps tempting to conjecture that the threshold for the appearance of long paths and cycles might also be controlled by the expansion ratio 
of the host graph, with a lower bound on their size being determined by the scale on which this expansion holds. Perhaps somewhat surprisingly, as noted in~\cite[Remark 3.3]{DEKK23+}, such a result cannot hold, even before percolation. Indeed, if we take our host graph to be a very unbalanced bipartite graph, for example, $G=K_{d,d^{10}}$, then whilst the graph has an expansion factor of at least $d$ for every set of size up to $d^{9}$, there is no path in $G$ of length longer than $2d$, since each path must have half of its vertices in the smaller partition class.

This suggests that one should perhaps consider the \emph{vertex-expansion} of the host graph, instead of the \textit{edge-expansion}. Indeed, our first result shows the typical emergence of long paths after percolation under the assumption of vertex-expansion at a fixed scale.
\begin{theorem}\label{t:longpath}
Let $k,d \in \mathbb{N}$, and let $G$ be a graph on at least $k$ vertices such that every $S\subseteq V(G)$ with $|S|= k$ satisfies $|N(S)| \geq kd$. Let $\varepsilon > 0$ be a sufficiently small constant and let $p = \frac{1+\varepsilon}{d}$. Then $G$ contains a path of length at least $\frac{\varepsilon^2 kd}{10}$ with probability at least $1- \exp\left( -\Omega_{\varepsilon}\left(k d \right)\right)$.
\end{theorem}
Note that for $k=1$ and $G=K_{d+1}$, the above recovers the result of Ajtai, Koml\'{o}s, and Szemer\'{e}di~\cite{AKS81a}. In fact, when $G$ is a clique of size $k(d+1)$, the above can be seen to be tight up to a multiplicative factor in the bound of the length of the path.

However, as is the case with arbitrary host graphs of large minimum degree, it is not immediately obvious how this result can be strengthened to find a long cycle. In the deterministic setting, that is, when considering a graph $G$ on at least $k$ vertices such that every $S\subseteq V(G)$ with $|S|=k$ satisfies $|N(S)|\ge kd$, it is known that $G$ contains a cycle of length $\Omega(kd)$ \cite{K19a}. In the specific case where $kd$ is linear in $n$, our main result shows that this holds after percolation as well.
\begin{theorem}\label{t:longcycle}
Let $0<c<1$, $d \in \mathbb{N}$, and let $G$ be a graph on $n$ vertices such that every $S\subseteq V(G)$ with $|S|= \frac{cn}{d}$ satisfies $|N(S)| \geq d|S|$. Let $\varepsilon\coloneqq \varepsilon(c) >0$ be a sufficiently small constant and let $p = \frac{1+\varepsilon}{d}$. Then $G_p$ contains a cycle of order $\Omega\left(\varepsilon^2c^2n\right)$ with probability at least $1-\exp\left(-\Omega_{\varepsilon,c}\left(n\right)\right)$.
\end{theorem}

As an application, we consider ``optimal'' $K_{s,t}$-free graphs, that is, graphs not containing a copy of $K_{s,t}$ (for $2\le s \le t$) and having maximal density. Tightly related to the Zarankiewicz problem, it is known that the maximal number of edges in an $n$-vertex $K_{s,t}$-free graph is $O(n^{2-\frac{1}{s}})$, and for $t$ sufficiently larger than $s$ (as well as for some fixed small values of $s$ and $t$), this has been established as the correct order of magnitude (see, for example, \cite{ARS99} and references therein). When a graph $G$ is an optimal $K_{s,t}$-free graph, it has sufficiently good vertex expansion properties \cite{KS09}, which allows one to derive that it has a cycle of length linear in $|V(G)|$ (utilising the connections between local expansion and long cycles as in \cite{K19a}). 

Let us show that in fact one can apply Theorem \ref{t:longcycle} to optimal $K_{s,t}$-free graphs. Indeed, let $\alpha$ be a positive constant, and let $G=(V,E)$ be a $K_{s,t}$-free graph with $|E|\ge \alpha n^{2-\frac{1}{s}}$. By a standard argument, there is a subgraph $G_0\subseteq G$ with minimum degree $\delta(G_0)\ge \alpha n^{1-\frac{1}{s}}$, and we note that $|V(G_0)|=\Theta(n)$ (indeed, otherwise $G_0$ will be too dense, contradicting it being $K_{s,t}$-free).
Recalling that any $K_{s,t}$-free graph $H$ has $O\left(|V(H)|^{2-\frac{1}{s}}\right)$ edges, for every $X\subseteq V(G_0)$ with $|X|=n^{1/s}$, we have $e(X,V\setminus X)\ge \frac{\alpha}{2}n=\frac{\alpha n^{1-1/s}}{2}|X|$. Hence, by~\cite[Lemma 7.2]{KS09}, for every $X\subseteq V(G_0)$ with $|X|=n^{1/s}$ we have $|N(X)|\ge\frac{\alpha n^{1-1/s}}{2t}|X|=\frac{\alpha}{2t}\cdot n$. We may thus apply Theorem~\ref{t:longcycle}, and obtain the following corollary:
\begin{corollary}
Let $\alpha$ be a positive constant, and let $s,t\in \mathbb{N}$ with $2\le s \le t$. Let $G$ be a $K_{s,t}$-free graph on $n$ vertices satisfying $|E(G)|\ge \alpha n^{2-\frac{1}{s}}$. Then, there exists a sufficiently large constant $K\coloneqq K(\alpha, t)$ such that for any $p\ge \frac{K}{n^{1-1/s}}$, with high probability $G_p$ contains a cycle of length $\Omega(n)$.
\end{corollary}


The paper is structured as follows. In Section \ref{s:prelim} we collect several definitions and results which will be useful throughout the paper. In Section \ref{s: path} we prove Theorem \ref{t:longpath}, and in Section \ref{s: cycle} we prove Theorem \ref{t:longcycle}. Finally, in Section \ref{s: discussion} we discuss our results and present avenues for future research.

\section{Preliminaries}\label{s:prelim}

\subsection{Notation}\label{s: notation}
Given subsets $A,B\subseteq V(G)$ with $A\cap B=\varnothing$ and a subgraph $H\subseteq G$, we denote by $e_H(A,B)$ the number of edges in $H$ with one endpoint in $A$ and the other endpoint in $B$. When $H=G$, we may omit the subscript. Furthermore, we denote by $N(A)$ the external neighbourhood of $A$ in $G$, that is, $N(A)\coloneqq \{v \in V(G)\setminus A \colon \exists u\in A, \{v,u\}\in E(G)\}$. For the sake of clarity we omit all rounding signs.

Given $I \subseteq \mathbb{N}$, we say that a graph $G$ is an $(I, d)$-expander if for every set $X$ with $|X| \in I$, it holds that $|N(X)| \geq d|X|$. If $G$ is a $(\{k\}, d)$-expander we may also simply say $G$ is a $(k, d)$-expander. We refer the reader to the surveys~\cite{HLW06,K19} for a comprehensive study of expander graphs and their applications.

\subsection{Depth-First Search (DFS)}\label{s:DFS}
The Depth-First Search (DFS) algorithm explores the components of a graph using a ``last-in-first-out'' discipline. The algorithm receives as input a graph $G=(V,E)$ and an ordering $\sigma$ of the vertex set $V$. During the algorithm, we maintain three sets of vertices: $W$, the set of \emph{processed} vertices; $A$, the set of \emph{active} vertices, which we treat as a stack; and $U$, the set of \emph{unvisited} vertices.

We initialise $W=A=\emptyset$ and $U=V$.
At each step, if the stack $A$ is non-empty, then we consider the most recently added $a\in A$, and query the edges from $a$ to $U$, according to the order $\sigma$, until a vertex $u\in U$ is discovered with $(a,u) \in E(G)$. If no vertex is found, $a$ is moved from the stack $A$ to $W$. Otherwise, we add the newly discovered vertex $u$ to the top of the stack $A$. If the stack $A$ is empty, we move the first vertex in $U$, according to $\sigma$, into $A$.

We note some elementary facts about this process:
\begin{enumerate}[(P\arabic*)]
\item\label{i:a} at each stage of the algorithm the stack $A$ spans a path in $G$; and,
\item\label{i:b} at each stage of the algorithm there are no edges in $G$ between $U$ and $W$.
\end{enumerate}

To analyse the DFS algorithm on a percolated graph $G_p$ we will utilise the \emph{principle of deferred decisions}. That is, we will sample a sequence $(X_i)_{i=1}^{|E(G)|}$ of i.i.d Bernoulli$(p)$ random variables, which we can think of as representing a positive (with probability $p$) or negative (with probability $1-p$) answer to a query in the algorithm. When the algorithm makes the $i$-th query (about an edge of $E(G)$), it receives a positive answer if and only if $X_i=1$. Let us denote by $W(i), A(i), U(i)$ the sets of vertices in $W,A,U$, respectively, after the $i$-th query. To complete the exploration of the graph, once 
$U$ is empty we make the algorithm query all the remaining edges in $G$ not queried before. It is then clear that the obtained graph has the same distribution as $G_p$.

\subsection{Concentration inequalities}\label{s:concentration}
We will make use of a typical Chernoff-type tail bound on the binomial distribution (see, for example,~\cite[Appendix A]{AS16}).
\begin{lemma}\label{l:Chernoff}
Let $n \in \mathbb{N}$, let $p \in [0,1]$ and let $X \sim \text{Bin}(n,p)$. Then, for every $0\le t \leq \frac{np}{2}$,
\[
\mathbb{P}\left[\left|X -np \right| > t\right] < 2 \exp\left(-\frac{t^2}{3np} \right).
\]
\end{lemma}

\section{Long paths: Proof of Theorem \ref{t:longpath}}\label{s: path} 
As discussed in the introduction, an assumption on edge-expansion alone does not suffice to prove Theorem \ref{t:longpath}. Nevertheless, the proof here is heavily inspired by that of~\cite{KS13}. 

We explore $G_p$ via the DFS algorithm described in \Cref{s:DFS} until it performs $N_1 = \frac{\varepsilon kd^2}{2(1+\varepsilon)}$ queries (indeed, $|E(G)|\ge kd^2/2>N_1$). The expected number of positive queries by this point is $\frac{\varepsilon k d}{2}$, and hence by Lemma \ref{l:Chernoff} there have been at least $\left(1-\frac{\varepsilon}{5}\right)\frac{\varepsilon kd}{2}$ positive queries with probability at least $1-\exp\left( -\Omega\left(\varepsilon^3 k d \right)\right)$. We will assume in what follows that this occurs deterministically.

Since either $U(N_1) = \emptyset$ (and $|V(G)| \geq kd$) or each positive query corresponds to a vertex which moves from $U$ to $A$ (which may later move to $W$), it follows that $|A(N_1) \cup W(N_1)| \geq  \left(1-\frac{\varepsilon}{5}\right)\frac{\varepsilon kd}{2}$. Let $N_2$ be the first time such that 
\begin{equation}\label{e:size}
    |A(N_2) \cup W(N_2)| = \left(1-\frac{\varepsilon}{5}\right)\frac{\varepsilon kd}{2},
\end{equation}
noting that $N_2 \leq N_1$. Suppose towards a contradiction that $G_p$ does not contain a path of length $\frac{\varepsilon^2 k d}{10}$. By \ref{i:a},  $|A(N_2)| \le \frac{\varepsilon^2 kd}{10}$, and therefore by \eqref{e:size}, $|W(N_2)| \geq \left(1-\frac{2\varepsilon}{5}\right)\frac{\varepsilon kd}{2}$. Let us fix disjoint sets $W_1, W_2, \ldots, W_r\subseteq W(N_2)$ of size $k$, where $r \geq \left(1-\frac{2\varepsilon}{5}\right)\frac{\varepsilon d}{2}$. Since $G$ is a $(k,d)$-expander, $|N(W_i)| \geq kd$ for each $i \in [r]$, and so by \eqref{e:size}
\[
    e(W_i,U(N_2)) \geq |N(W_i)|-|A(N_2)\cup W(N_2)|\ge kd - |A(N_2) \cup W(N_2)| = \left(1-\frac{\varepsilon}{2}+\frac{\varepsilon^2}{10}\right)kd.
\]
Since the sets $W_i$ are disjoint, it follows that 
\[
e(W(N_2),U(N_2)) \geq r \left(1-\frac{\varepsilon}{2}+\frac{\varepsilon^2}{10}\right)kd\geq \left(1-\frac{\varepsilon}{2}\right)\left(1-\frac{2\varepsilon}{5}\right) \frac{\varepsilon kd^2}{2}
\]
However, by \ref{i:b} every edge between $U(N_2)$ and $W(N_2)$ has already been queried and hence
\[
\left(1-\frac{\varepsilon}{2}\right)\left(1-\frac{2\varepsilon}{5}\right) \frac{\varepsilon kd^2}{2} \leq e(W(N_2),U(N_2)) \leq N_2 \leq N_1 = \frac{\varepsilon kd^2}{2(1+\varepsilon)},
\]
or equivalently $(1+\varepsilon)\left(1-\frac{\varepsilon}{2}\right)\left(1-\frac{2\varepsilon}{5}\right) \leq 1,$ which is a contradiction for $\varepsilon$ sufficiently small. \qed
\section{Long cycles: Proof of Theorem \ref{t:longcycle}}\label{s: cycle}

The proof will proceed via a two-round exposure process, using $p_1= \frac{1+\frac{\varepsilon}{2}}{d}$ and $p_2 = \frac{p-p_1}{1-p_1} \geq \frac{\varepsilon}{2d}$, noting that $G_p$ has the same distribution as $G_{p_1}\cup G_{p_2}$, since $1-p=(1-p_1)(1-p_2)$.
  
Let $\alpha \coloneqq \frac{c(\varepsilon/2)^2}{10}=\frac{c\varepsilon^2}{40}$. By \Cref{t:longpath}, $G_{p_1}$ contains a path $P$ of length $\alpha n$ with probability at least $1- \exp\left( -\Omega_{\varepsilon}\left(n\right)\right)$. We continue assuming that this holds deterministically. Let us fix $r\coloneqq\frac{\alpha d}{c}$ disjoint subpaths $P_1\cup \ldots \cup P_r$ of $P$, which we call \emph{blocks}, each of length $\frac{cn}{d}$. We may assume without loss of generality that these blocks are ordered according to the order they appear in the path $P$.

For each $v \in V(G) \setminus V(P)$, let $B_v = \{ i \in [r] \colon N(v) \cap V(P_i) \neq \varnothing\}$ be set of indices of blocks to which $v$ is adjacent (in $G$). We call $v$ \emph{good} if $b_v \coloneqq |B_v| \geq \frac{\alpha d}{2}$, and denote the number of good vertices by $g$.

By our assumption on $G$, every block $P_i$ has at least $d \cdot |P_i| - |P| = \left(c - \alpha\right)n$ neighbours outside $P$, and in total there are  $r = \frac{\alpha d}{c}$ blocks. Hence, by a double-counting argument 
\[ \left(c - \alpha\right)n\cdot \frac{\alpha d}{c} \le \sum_{v \in V\setminus P} b_v \le     g\cdot \frac{\alpha d}{c} + n \cdot \frac{\alpha d}{2}.\]
In particular, $g \geq \left(\frac{c}{2} - \alpha\right)n \geq \frac{cn}{3}$, since $\varepsilon$ is sufficiently small.

We say a good vertex is \emph{successful} if there is at least one edge of $G_{p_2}$ between $v$ and a vertex $v_f$ of the first $\frac{b_v}{3}$ blocks of $B_v$ \textit{and} one edge between $v$ and a vertex $v_\ell$ of the last $\frac{b_v}{3}$ blocks of $B_v$. Note that a successful vertex $v$ lies in a cycle in $G_{p_1}\cup G_{p_2}$ of length at least $\frac{\alpha c n}{6}$. Indeed, by assumption there are $i$ and $j$ with $v_f \in P_i$, $v_\ell \in P_j$, $\{\{v,v_f\},\{v,v_\ell\}\} \subseteq E(G_{p_2})$ and $|i-j| - 1 \geq \frac{b_v}{3} \geq  \frac{\alpha d}{6}$. Since $v \not\in V(P)$, the vertex $v$ together with (part of) $P$ forms a cycle in $G_{p_1}\cup G_{p_2}$ of length at least $\frac{\alpha d}{6}\cdot \frac{cn}{d}=\frac{\alpha cn}{6}=\Omega(\varepsilon^2c^2n)$.

Each good vertex $v$ satisfies $b_v \geq \frac{\alpha d}{2}$ by definition, and hence
\[
  p_{\text{s}} \coloneqq \mathbb{P}[v\text{ is successful}] \geq \left(1 - (1-p_2)^{\frac{b_v}{3}}\right)^2 
  \geq \left(1 - \exp\left(-\frac{\varepsilon \alpha}{12}\right)\right)^2= \Omega(c^2\varepsilon^6).
\]
There are $g\ge \frac{cn}{3}$ good vertices, and each good vertex is successful independently of the others. Therefore, the probability that no good vertex is successful is $(1-p_{\text{s}})^{g} = \exp\left(-\Omega_{\varepsilon,c}\left(n\right)\right)$.
We conclude that with probability at least $1-\exp\left( -\Omega_{\varepsilon}\left(n\right)\right)-\exp\left(-\Omega_{\varepsilon,c}\left(n\right)\right)=1-\exp\left(-\Omega_{\varepsilon,c}\left(n\right)\right)$, there is a cycle of length at least $\Omega(\varepsilon^2c^2n)$ in $G_{p_1}\cup G_{p_2}$.\qed

\section{Discussion}\label{s: discussion}
Whilst Theorem \ref{t:longcycle} requires that the graph expands at an almost optimal scale, that is $d$-vertex-expansion for sets of size $\Omega(n/d)$, we conjecture that the conclusion should hold for any $(k,d)$-expander.
\begin{conjecture}\label{conjecture}
Let $k,d \in \mathbb{N}$ and let $G$ be a graph on at least $k$ vertices such that every $S\subseteq V(G)$ with $|S|= k$ satisfies $|N(S)| \geq kd$. Let $\varepsilon > 0$ be sufficiently small and let $p = \frac{1+\varepsilon}{d}$. Then with probability tending to one as $k$ tends to infinity, $G_p$ contains a cycle of length $\Omega_{\varepsilon}(kd)$.
\end{conjecture}
Note that Conjecture \ref{conjecture} would be a strengthening of Theorem \ref{t:longpath}, where we show that typically $G_p$ contains a path of length $\Omega_{\varepsilon}(kd)$.

As noted in the introduction, while assuming that the minimum degree is at least $d$ suffices to show that $G_p$ contains a cycle of length linear in $d$ \cite{KS14}, an assumption on the edge-expansion alone does not suffice to prove Theorem \ref{t:longpath}. Such an assumption is tightly related to an assumption on $\lambda_2$, the second largest eigenvalue of the adjacency matrix of the graph. In contrast, it is known that an appropriate assumption on $\lambda=\min\{|\lambda_2|,|\lambda_{|V(G)|}|\}$ does suffice. Indeed, for $(n,d,\lambda)$-graphs $G$ with growing degree $d$, when $p=\frac{1+\varepsilon}{d}$ and $\lambda\le \varepsilon^4d$, Diskin and Krivelevich~\cite{DK21B} showed that \textbf{whp} $G_p$ contains a cycle of length linear in $n$. 

Considering constant $d$, Alon and Bachmat~\cite{AB06} showed that given a $(d+1)$-regular graph $G$ on $n$ vertices, for any fixed $\varepsilon >0$, if $p=\frac{1+\varepsilon}{d}$, then with probability at least $1-o_n(1)$ the random subgraph $G_p$ contains a cycle. Theorems \ref{t:longpath} and \ref{t:longcycle} can be seen as a quantitative strengthening of this result, showing what additional requirements are sufficient to show the typical existence of paths or cycles of some given length. 



\section*{Acknowledgement} The first and third authors were supported in part by the Austrian Science Fund (FWF) [10.55776/\text{P36131}]. Part of this work was done while the first and the third authors were visiting Tel Aviv University, and they would like to thank the university for its hospitality. For the purpose of open access, the authors have applied a CC-BY public copyright licence to any Author Accepted Manuscript version arising from this submission.

\bibliographystyle{abbrv}
\bibliography{bigbib}

\end{document}